\documentclass{article}
\usepackage{url} 
\usepackage{soul}

\usepackage{caption}
\usepackage{subcaption}
\usepackage{amsmath}
\usepackage{amsthm}
\usepackage{amssymb}
\usepackage{amsfonts}
\usepackage{mathtools}
\usepackage{thmtools}
\usepackage{graphicx}

\usepackage{hyperref}
\usepackage[noabbrev,nameinlink,capitalise]{cleveref}

\newcommand{\define}[1]{{\bf \boldmath{#1}}}

\usepackage[backend=biber]{biblatex}
\title{SuPerPoV: Score and evolution of the stratospheric polar vortex via persistent homology}

\author{Jake Cordes, Barbara Giunti, Zheng Wu}
\date{}

\begin{document}

\maketitle

\begin{abstract}
Classifying the stratospheric polar vortex provides predictability for surface weather on extended-range timescales. However, providing a scientifically sound classification is challenging: all the definitions proposed in over 60 years of study depend on empirically chosen parameters and yield different results when one of them changes. Moreover, as they are based on static thresholds, it is not straightforward to use them to study the spatiotemporal evolution of the vortex. Here, we introduce SuPerPoV, a score system that computes displacement and split ratios of the polar vortex using tools from topological data analysis, thus providing a sound classification of the polar vortex. The scores are computed by adapting superlevel set persistence and comparing prominent features. Our definition is entirely threshold-free and implemented open source. The scores generally recovers previous definitions and are output for a user-defined number of days, thus showing the evolution of the event. SuPerPoV offers a paradigm shift in the study of the polar vortex, hopefully bringing a deeper understanding of the polar vortex and related extreme events, such as sudden stratospheric warmings.
\end{abstract}

\section{Introduction} 
Sudden stratospheric warmings (SSWs) were first observed as ``explosive warmings in the stratosphere" in 1952, accompanied by deceleration of the \define{polar vortex}, which consists of strong circumpolar westerly winds in the polar stratosphere that form in fall and decay in spring \cite{baldwin2021sudden}. 
SSWs can have a huge impact on the tropospheric circulation and surface weather for up to 2 months after the onset of events \cite{baldwin2001stratospheric,kidston2015stratospheric}, and thus are considered as an important source of surface predictability on subseasonal to seasonal (S2S) timescales \cite{karpechko2017predictability}. 
Since SSWs occur almost exclusively in the Northern Hemisphere (NH), we focus our analysis on the stratospheric polar vortex in the NH.
The deceleration and distortion of the polar vortex are caused by the interaction with the upward-propagating planetary waves \cite{matsuno1971dynamical}. 
Based on the polar vortex geometry at the onset of SSWs, they can generally be classified into \define{displacement events}, when the vortex is shifted off the pole, and \define{split events}, when the vortex splits into two smaller vortices \cite{charlton2007new}. 
In addition, SSWs can also be categorized into minor and major events based on the magnitude of the polar vortex disruption \cite{butler_defining_2015}.
If there is no SSW event, then the polar vortex remains unchanged and is \define{normal}. 
Thus, classifying the structure and the degree of disruption of the stratospheric polar vortex is key to understanding and identifying SSWs and their surface impact.
\newline
Despite decades of research into SSWs, various, discordant definitions have been proposed based on differing criteria and empirical fixed thresholds that remain somewhat subjective and arbitrary. 
Traditional frameworks, such as the widely adopted metric of zonal-mean zonal wind reversal at 10 hPa, often impose discrete cut-off points on a continuous dynamical spectrum of vortex states. 
Meanwhile, the prior definitions are also sensitive to predefined thresholds, and small adjustments to these parameters can lead to different classification results, as shown in \cite{butler_defining_2015}.
Moreover, these definitions often require smooth inputs, and the choices of the parameters are hidden, making fine-tuning challenging. Therefore, ongoing efforts are required to refine a more robust and objective definition \cite{butler_defining_2015}. 
On the other hand, intuitively, it should not be particularly difficult to differentiate between a split vortex (there are two ``cones''), a displaced one (the ``cone's tip'' is not on the North Pole), and a normal one (the ``cone's tip'' is right above the North Pole). 
In particular, when we think about the polar vortex as a \emph{shape}, classifying it seems straightforward.
\newline
Applied topology is a field at the intersection of topology, computer science, and data analysis that uses tools from topology to discover underlying features in the data, where topology is a branch of mathematics related to geometry that focuses on properties of shapes that remain unchanged under continuous transformations of the object.
Thus, it seems only logical to use tools from applied topology to study the stratospheric polar vortex and other atmospheric science phenomena. 
Indeed, there have already been successful uses of applied topology in this area \cite{licon-salaiz_2018,muszynski2019topological,tymochko_2020}, even if, to the best of our knowledge, none about polar vortices and SSWs.
In this work, we show how \emph{persistent homology}, one of the most popular tools of applied topology, provides a nuanced classification of polar vortices, distinguishing between split and displaced, and interpolating between all previous definitions from the literature.
\newline
Homology is the topological tool that identifies ``holes'' in the shapes: a $0$-hole is the blank space between the components of an object, a $1$-hole can be thought as the ``loop'' formed by a rubber band - it may not be a circle, but it is \emph{circular} -, a $2$-hole is the cavitity inside a football, and so on. 
The polar vortex will be visualized and analyzed by scanning the scalar field of the geopotential height (gph) over the NH using persistent homology, with different holes representing different vortex states. 
The output of the persistent homology is two scores (split and displacement) that can inform us of the persistence and transition of vortex states, allowing end-users to identify SSWs themselves.
The split and displacement scores are very straightforward to interpret, and their computation is objective since it does not rely on any parameter. 
Moreover, while persistent homology is used in their computation, there is no need to know what persistent homology is to interpret the scores. 
The algorithm computing the scores is very clean and short, and can be run with ease even by non-experts. 
All instructions are provided at \cite{github_repo}.
Therefore, our method can be used without any knowledge of applied topology. 
Since the output is in the form of scores, the users will need to decide at what scale one has a major SSW. 
However, this approach provides a more comprehensive view of vortex evolution and transitions, enabling the identification of early signals for SSW events, and can be tailored to the needs at hand. 
Moreover, this output format proves that our method interpolates all the previous definitions: using a 10 hPa gph, events that were categorized by many definitions in \cite{butler_defining_2015} have a large percentage of splitness or of displacement, while events that were categorized only by a few definitions in the literature have a lower value in both scores (but never zero).
Even if most of the analysis we conducted used ghp at 10 hPa, this is not a requirement. 
Indeed, we added a comparison with gph at 50 and 100 hPa to better understand how polar vortices evolve in the lower stratosphere. 
Last but not least, when running our analysis, the user can decide the span of days they want to obtain the score for, and thus can see the temporal evolution of the vortex on shorter and longer scales. 
We hope that this will aid the prediction of polar vortex states and SSWs and potentially shed light on their impact on surface weather.
\section{Data and Methodology}

\subsection{Data visualization}

\begin{figure}[!h]
    \centering
    \begin{subfigure}{0.32\textwidth}
        \centering
        \includegraphics[width=\linewidth]{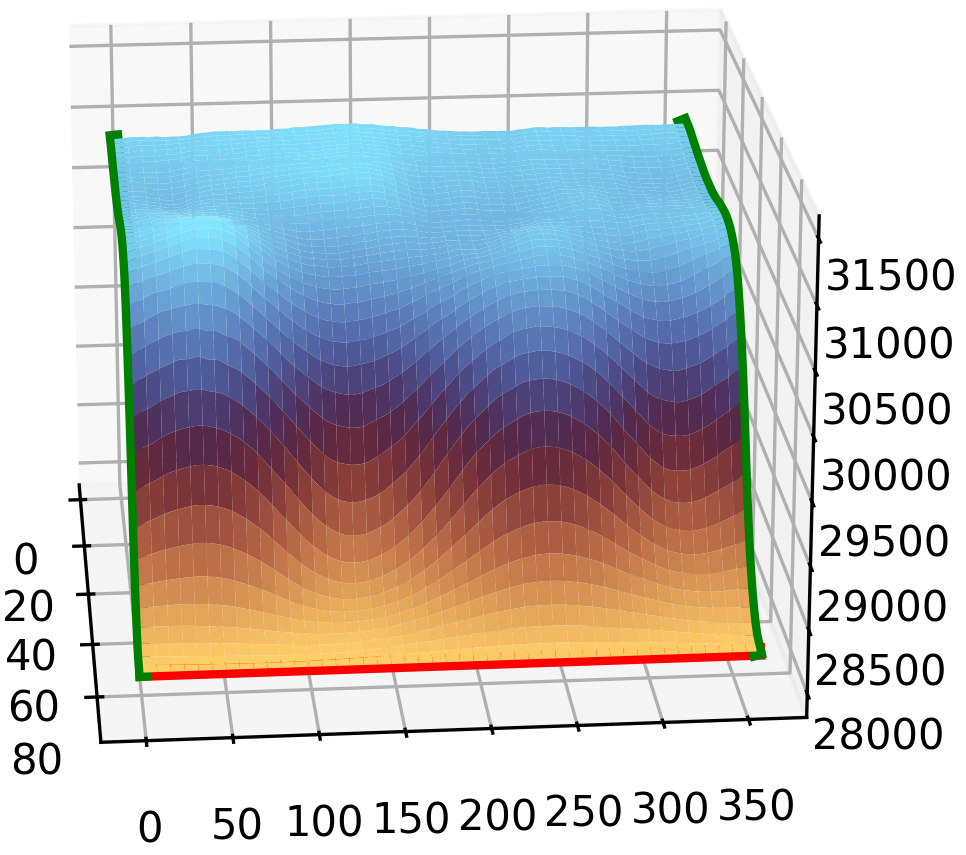} \\
        \includegraphics[width=\linewidth]{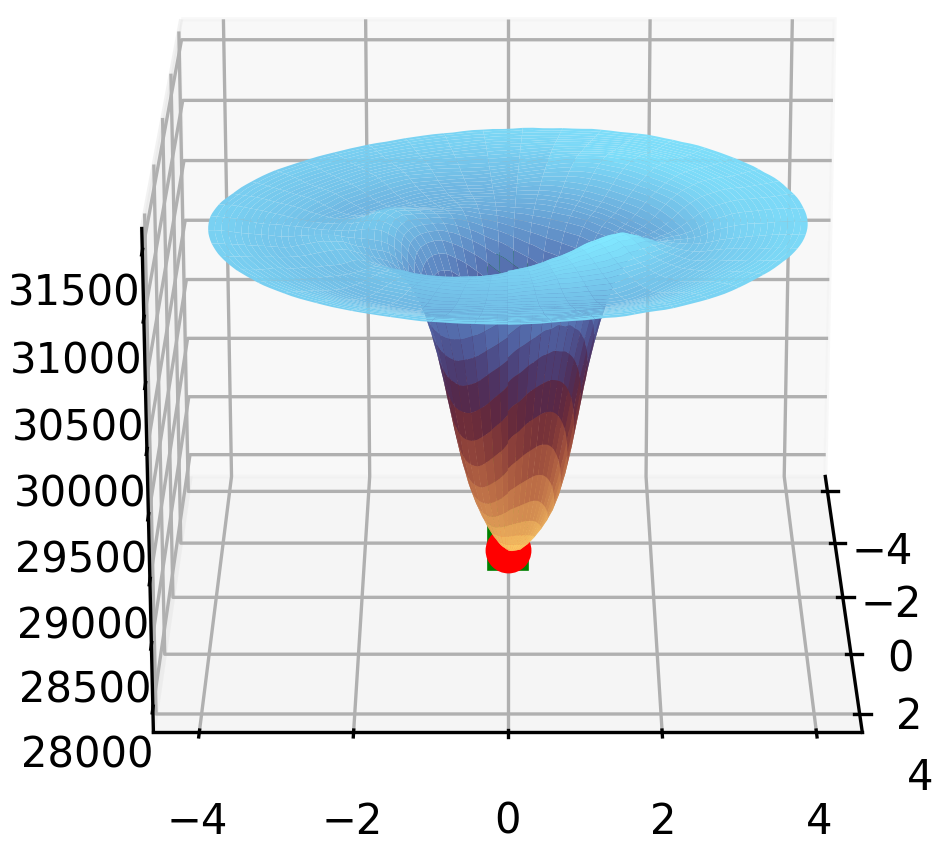}
        \caption{1961-02-08}
        \label{subfig:normal}
    \end{subfigure}
    \hfill
    \begin{subfigure}{0.32\textwidth}
        \centering
        \includegraphics[width=\linewidth]{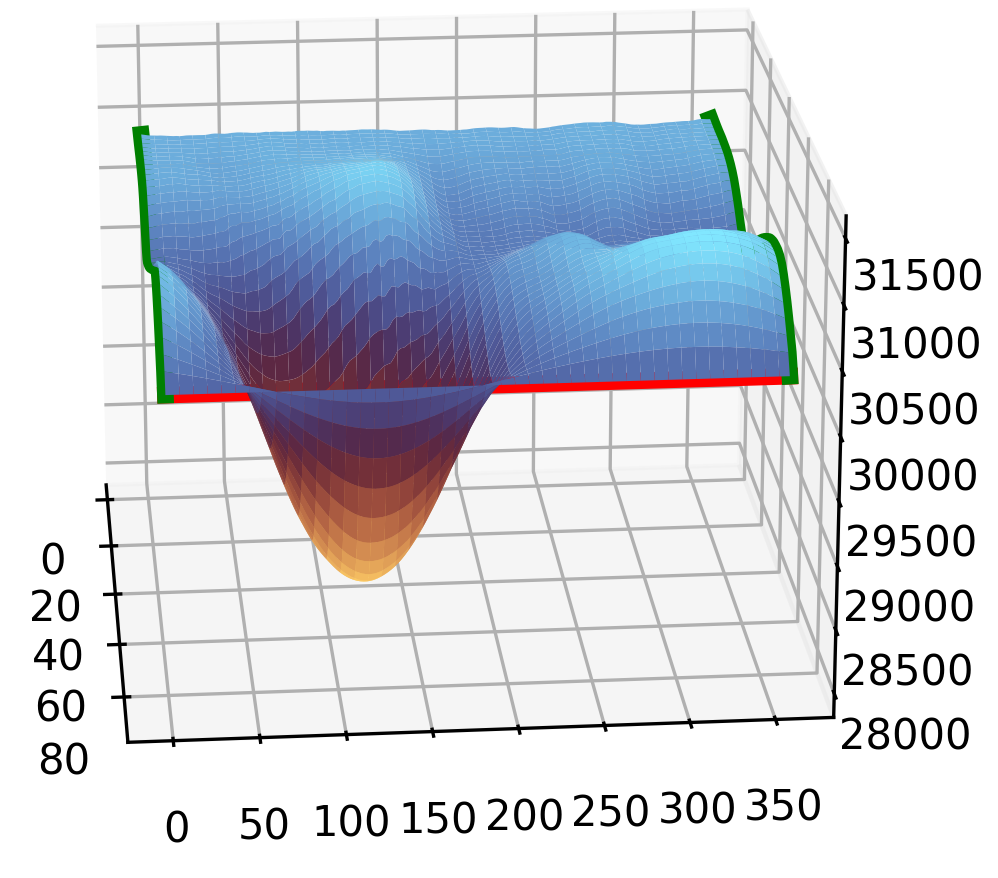} \\
        \includegraphics[width=\linewidth]{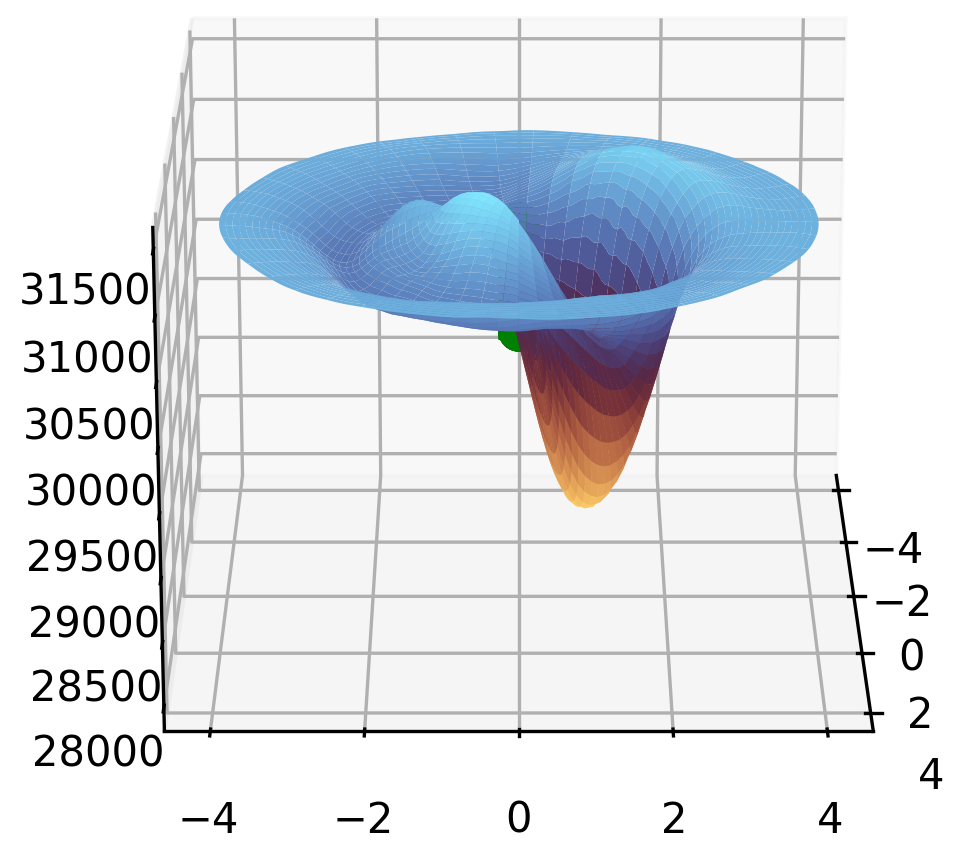}
        \caption{1965-12-23}
        \label{subfig:displaced}
    \end{subfigure}
    \hfill
    \begin{subfigure}{0.32\textwidth}
        \centering
        \includegraphics[width=\linewidth]{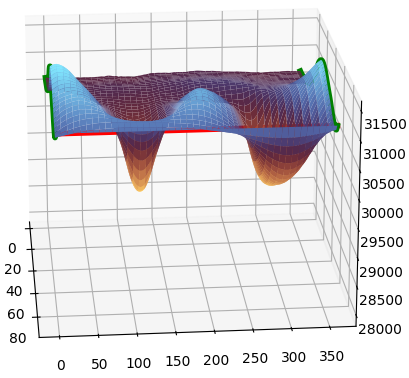} \\
        \includegraphics[width=\linewidth]{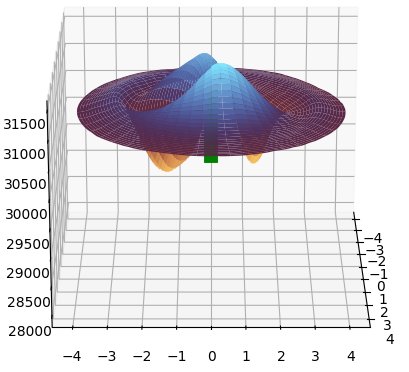}
        \caption{1989-02-21}
        \label{subfig:split}
    \end{subfigure}
    \caption{Grid (top row) and cylindrical (bottom row) plots for a normal vortex (a), a displaced one (b), and a split one (c). The red line/dot denotes the North Pole, and the green lines in the grid plot are identified in the cylindrical one.}
    \label{fig:plots}
\end{figure}

The data used in the study is the daily \textbf{geopotential height (gph)} over the NH at a fixed pressure in the stratosphere during the extended winter from November to April from 1959 to 2021 of ERA5 reanalysis from the European Center for Medium-Range Weather Forecasts (ECMWF) \cite{hersbach2020era5}. The gph at 10 hPa is the primary variable used to represent the overall shape of a vortex and analyzed in this study. Besides, we also use the gph at 50 hPa and 100 hPa to further examine the vortex states in the lower stratosphere and the feasibility of using them to identify major SSW events. 
We examine the data in two ways: 
\begin{itemize}
    \item \define{Grid Plots} are a visualization of the data plotted in the XYZ plane, and the points are of the form (lat, lon, gph), see Figure \ref{fig:plots}, top row.
    The North Pole is represented as a red line.   
    \item \define{Cylindrical Plots} are a visualization of the data where we convert the latitude-longitude into polar coordinates, and plot the height values as the third coordinate, see Figure \ref{fig:plots}, bottom row.
    
\end{itemize}
Several definitions in the literature used, together or instead of the gph, the \define{zonal-mean zonal wind speed} at 60$^\circ$ N and 10 hPa, to classify polar vortices. 
While our computation does not take it into account, we included into our visualization for comparison with other methods.

\subsection{Persistent Homology}

\begin{figure}[!h]
    \centering
    \begin{subfigure}{0.19\textwidth}
        \centering
        \includegraphics[width=\linewidth]{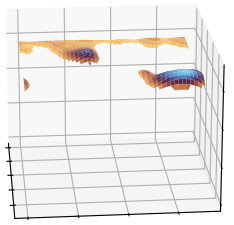}
        \includegraphics[width=\linewidth]{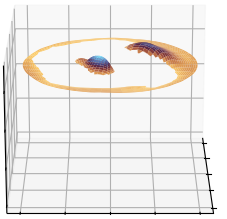}
        \caption{30830}
        \label{subfig:30830}
    \end{subfigure}
    \hfill
    \begin{subfigure}{0.19\textwidth}
        \centering
        \includegraphics[width=\linewidth]{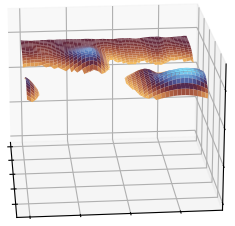}
        \includegraphics[width=\linewidth]{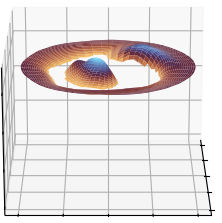}
        \caption{30600}
        \label{subfig:30600}
    \end{subfigure}
    \hfill
    \begin{subfigure}{0.19\textwidth}
        \centering
        \includegraphics[width=\linewidth]{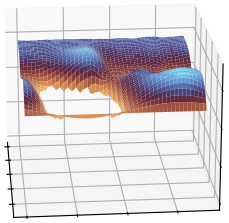}
        \includegraphics[width=\linewidth]{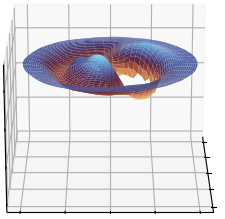}
        \caption{30150}
        \label{subfig:30150}
    \end{subfigure}
    \hfill
    \begin{subfigure}{0.19\textwidth}
        \centering
        \includegraphics[width=\linewidth]{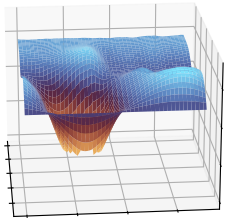}
        \includegraphics[width=\linewidth]{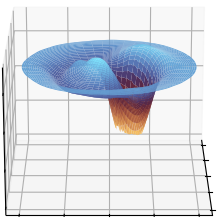}
        \caption{29000}
        \label{subfig:29000}
    \end{subfigure}
    \hfill
    \begin{subfigure}{0.19\textwidth}
        \centering
        \includegraphics[width=\linewidth]{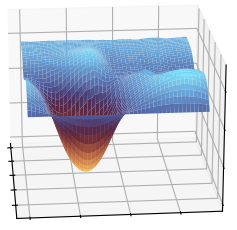}
        \includegraphics[width=\linewidth]{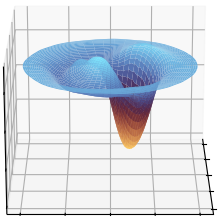}
        \caption{21000}
        \label{subfig:21000}
    \end{subfigure}
    \caption{A visualization of superlevel-set filtration on the grid (top row) and cylindrical (bottom row) plot for the date 1965-12-23. The value at the bottom of each column is the height we have reached going downward (so, everything above it is included).}%
    \label{fig:superlevelset_filt}
\end{figure}

We will provide a high-level intuition of persistent homology here, and refer the reader to standard textbooks \cite{Dey_Wang_2022,monster} for more details.
Recall that we will be computing persistent homology on the plots defined above (shown in Figure \ref{fig:plots}). 
In particular, the plots in circular coordinates let us visualize polar vortices as ``cones''. 
We can imagine slicing them horizontally at a certain height: if we find two $1$-holes, the polar vortex is split, while only one $1$-hole is not conclusive; it could be normal or displaced, or it could be that one of the two cones of a split polar vortex starts higher up. 
Therefore, ideally, we would like to scan the polar vortex downward and record when and where there are $1$-holes, to find where exactly the main $1$-hole is (above the North Pole or not) and to see if there are ever two $1$-holes. 
This is exactly what (superlevel set) persistent homology does: it computes the homology not just at one parameter (the slicing's height), but for all the gph's heights. 
The output is then not just the presence of holes, but also their evolutions. 
\newline
We focus on the $1$-holes, and record their lifespan: from the height at which a $1$-hole appears for the first time while scanning the geopotential scalar field downward, to the height at which that $1$-hole is filled (i.e., we reached the ``tip of the cone'').
These lifespans will give us the SuPerPoV split and displacement scores (see Section \ref{sec:superpov}).
\newline
Formally, to compute persistent homology, the first thing we need is a \textbf{filtration}.
A filtration on a (topological) space $X$ is a sequence of nested subspaces $\{X_i\}_{i=1}^{n}$ in $X$ such that $X_1 \subset X_2 \subset \dots \subset X_n \subseteq X$.
The filtration we use in our analysis is a \define{super-levelset filtration}, where we scan the surface (i.e., the space) downward.
We underline that the surface needs \textbf{not} to be smooth. 
Therefore, we can use any triangulation straightforwardly.
Such filtration can be thought of as the surface being completely submerged in coffee, and we examine the portions (i.e., the nested subspaces) that reveal themselves as we drain the coffee (see Figure \ref{fig:superlevelset_filt}).
For each subspace in the filtration, we compute its \define{$1$-homology} $H_1$, a tool from algebraic topology.
We adapt the mathematical language and say that the space at that height has ``an $H_1$ cycle'' or simply ``an $H_1$''.\footnote{This is technically incorrect; the precise formulation is ``a generator of $H_1$''.}
As we scan down a vortex, some $H_1$ cycles will appear at certain heights and disappear at lower heights.
The height at which an $H_1$ forms is called the \define{birth height} of the $H_1$, while the height at which an $H_1$ cycle disappears is called its \define{death height}\footnote{In applied topology literature, these are more commonly referred to as birth and death \emph{times}.}. 
Similarly, the \define{lifespan} of an $H_1$ (or \define{$H_1$-lifespan}), is given by the difference between the birth and the death heights of a given $H_1$. 
Since, in general, any plot will have multiple $H_1$, we can order them by their lifespans. 
In this work, we are interested in the two longest $H_1$-lifespans, but in the future, more fine-tuned analysis may also encompass shorter lifespans to obtain a more refined classification.
Figure~\ref{fig:superlevelset_filt} shows an example of applying super-levelset filtration to the 10 hPa gph plots on 1965-12-23. 
There, one notices that an $H_1$ is formed in the grid plot at height 30150, while an $H_1$ was formed in the cylindrical plot at height 30830.
In both cases, the $H_1$ persists until the very end of the filtration, where the hole then gets filled up. 
Note that there could be multiple $H_1$ in each plot, especially for split vortices. 
We pick the longer-lasting $H_1$ and use it to define the height of a vortex, both in the grid and in the cylindrical plot, which we denote by $h_{\mathrm{grid}}$ and $h_{\mathrm{cyl}}$, respectively.
\subsection{SuPerPoV scores}\label{sec:superpov}
We now have all the ingredients needed to define our novel contributions, the SuPerPoV (SUperlevel PERsistence of POlar Vortices) scores for analyzing and identifying the vortex states and SSW events. 
The \define{split score} represents how split the vortex is, and the \define{displacement score} represents how displaced the vortex is.
If both scores are close to zero, then we will conclude that the vortex is in a normal state. 
On the other hand, the larger the scores are, the more intense the displacement and/or splitness are.
The split score is given by the ratio between the longest $H_1$-lifespan and the second longest $H_1$-lifespan. 
\begin{equation}
    \text{Split Score} = \dfrac{\text{Second longest } H_1\text{-lifespan}}{\text{Longest } H_1\text{-lifespan}} \, .
\end{equation}
Intuitively, if the vortex is not split, the second longest lifespan will be very small, while if the vortex is truly split, the two $H_1$-lifespans will have more or less the same size, i.e., the two $1$-holes have a lifespan of similar length (see Figure \ref{subfig:split}). 

The displacement is also measured with a score, between the longest lifespan computed in two different coordinate systems: the original grid and the circular plot given by the input in polar coordinates
\begin{equation}
    \text{Displacement Score} = \dfrac{h_\mathrm{grid}}{h_\mathrm{cyl}} \, .
\end{equation}
A normal polar vortex appears as a tilted sheet of paper in the grid coordinate, and as a sharp cone in the circular plot (see Figure \ref{subfig:normal}). 
Thus, the lifespan of the longest-lasting $1$-hole will be close to zero in the former case and far from zero in the latter, and their ratio is very small. 
On the other hand, the lengths of the longest-lasting $1$-holes in the two plots are very similar if the vortex is displaced, and their ratio is very close to or above $1$ (see Figure \ref{subfig:displaced}).

\section{Usage and results}

The code computing the SuPerPoV scores is freely available at \cite{github_repo}. 
In it, we call modules from GUDHI \cite{gudhi} to compute persistent homology. 
While all the data we used is available at \cite{our_zenodo}, we remark that the code can be used on any gph grid, for any (collection of) day(s), and at different pressure levels, and thus can be applied beyond the scope of this paper in different applications that may not be related to the stratospheric polar vortex. 
The software asks the user for which day they are interested in classifying and analyzing the polar vortex state, as well as for how many days before and after they want to obtain. 
If one is interested only in one specific day, they can put 0 for the number of days before and after, but the software supports the computation for whole winters.
Note that, in the case of many days, some of the day labels are skipped for readability, but all scores are plotted.

\begin{figure}[!h]
    \centering
    \includegraphics[scale = 0.37]{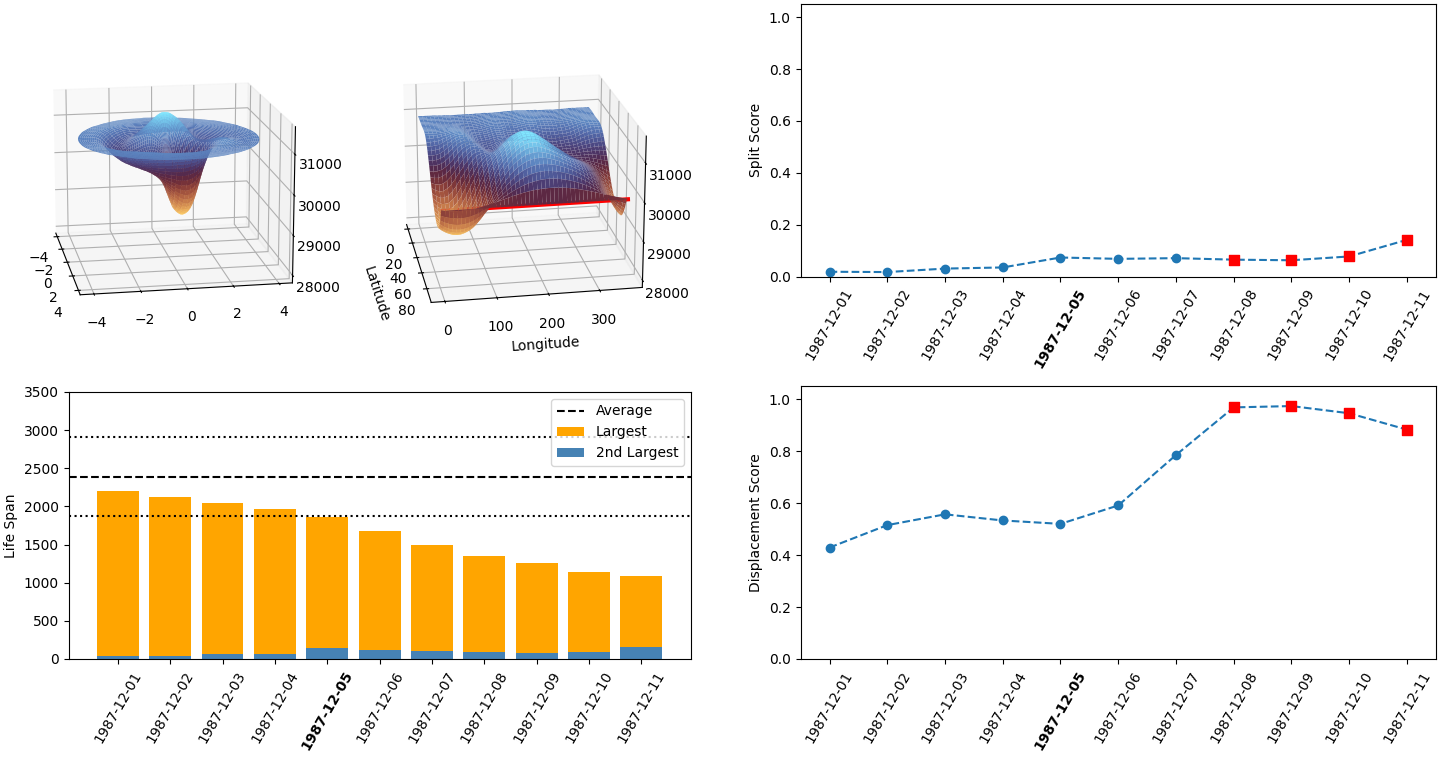}
    \caption{SuPerPoV output of the 1987-12-5, with 5 days prior and 6 days after. Grid and cylindrical plots of the day in the top left corner, lifespans of the two longest $H_1$ for each day in the bottom left, split and displacement scores on the right. In the latter two, red points denote negative zonal-mean wind speed. In all plots, the chosen input day is highlighted in bold.}
    \label{fig:superpov_output}
\end{figure}

\subsection{Output visualization and interpretation}

To demonstrate the SuPerPov output and its interpretation, we use one case study as an example, namely the gph at 10hPa on the 5th of December 1987, with 4 days prior and 6 days after. 
The output is depicted in Figure~\ref{fig:superpov_output}. 
In the labels, the input day is in bold.
On the top left, we plotted the grid and cylindrical plots of the input day (but not of the previous/following days) to aid understanding. 
We can see the displaced structure of the polar vortex.  
In the bottom left corner, we plotted the lifespans of the longest and second longest $H_1$-lifespans, together with a dashed line for the average height of normal vortices and two dotted lines representing one standard deviation above and below the average, respectively\footnote{For plotting these lines, we took the lifespans of the polar vortex for all days that had less than 0.1 and its split score was less than 0.05. While these values were somewhat arbitrarily chosen, they are used only for visualization purposes and not for computing the SuPerPov scores.}. 
We notice that, in general, both displaced and split polar vortices were below the average $H_1$-lifespan, suggesting that whenever the polar vortex is perturbed, it becomes shallower, i.e., less deep.
This plot also allows us to gauge when a polar vortex is particularly shallow, without being split, and remains centered on the North Pole.  
Since this tends to happen in late March and April, with the usual seasonal weakening of the polar vortex, our output reflects the seasonal transitions of the polar vortex.
The split and displacement scores are plotted one above the other on the right. 
We can see the displacement score slowly rising from 0.4 (a moderately displaced event), oscillating slightly, until Dec 6th, and then jumping quickly to almost 1, when the polar vortex is intensely displaced. 
On the other hand, the split score remains below 0.3 during the whole period, but it is decidely increasing toward Dec 12th, suggesting that the vortex is about to split. 
Moreover, the high displacement score is consistent with the low lifespan values plotted on the left.
In these plots, the red squares denote days in which the zonal-mean zonal wind speed at latitude 60$^\circ$ N and 10 hPa blows eastward\footnote{because this was our input, our code would run equally using as an input a different wind speed.}. 
We added this feature to facilitate the interpretation between high displacement score and the negative zonal-mean windspeed, but it is not needed in the score computation, and if this data is not available, it will not be plotted without error.
As Figure~\ref{fig:superpov_output} conveys, studying the evolution of the polar vortex is very natural with our method.
\begin{figure}[!h]
    \centering
    \includegraphics[scale = 0.4]{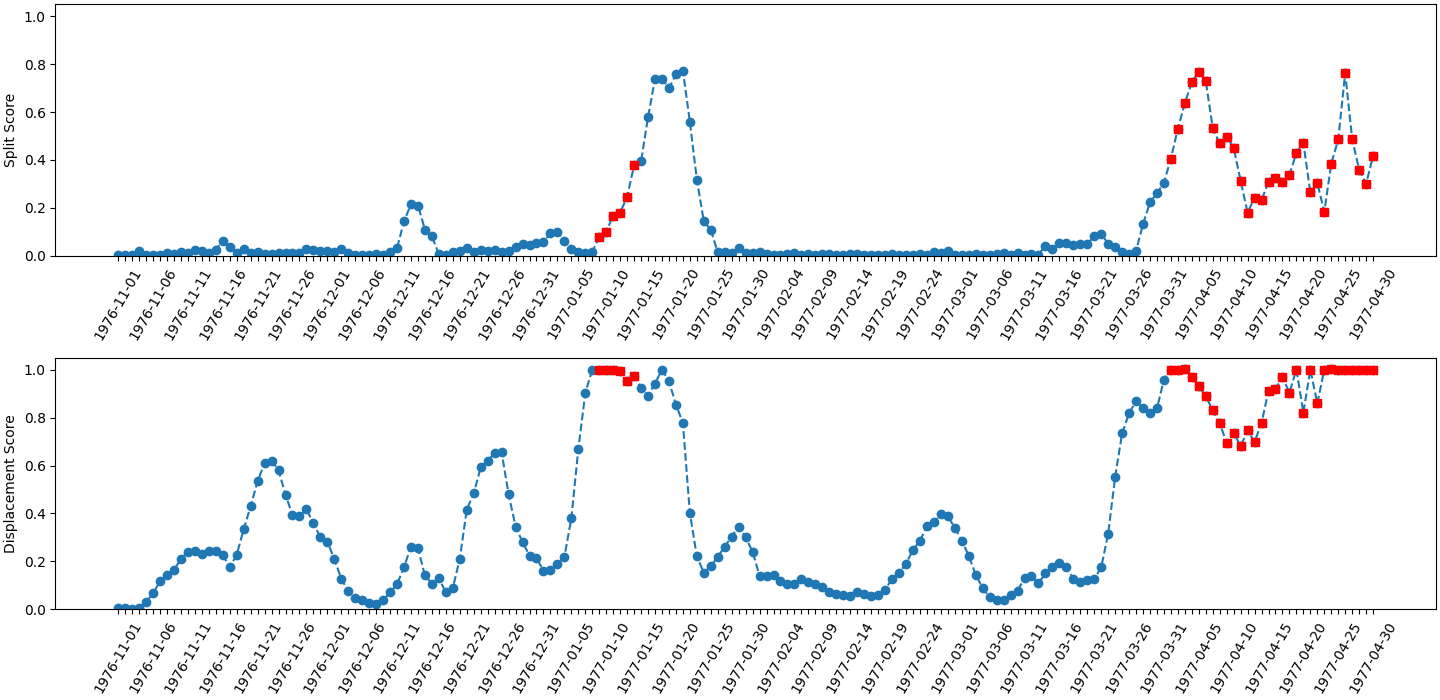}
    \caption{SuPerPoV scores of the whole winter 1976-1977. The red points represent a negative zonal-mean wind speed, to facilitate the comparison with other methods.}
    \label{fig:winter_76}
\end{figure}
To further showcase how our method provides a novel framework for studying the continuous behvior of polar vortex, in Figure~\ref{fig:winter_76} we plot the evolution of the SuPerPoV scores for the whole winter of 1976-77, which was particularly turbulent, with several displacement and split vortex states. 
While the decision of at which precise score value an event occurs is deferred to the case at hand, it is clear that some form of displacement happened around November 21st, December 21st, January 5th (this lasting for until the end of January), possibly around February 24th, and again at the end of March. 
A minor split event may have occurred on the 11th of December (resulting also in a small increase of the displacement score), while around January 15th and the end of March we see two clear split events happening.
Consistent with previous definitions, the last days in late March and April often see major activities, but this is most likely due to the seasonal transition.

\subsection{Comparison with previous methods}
Lastly, we compared the SuPerPoV scores with previous methods. 
To do so, we took all days classified as major SSW events by all previous definitions listed in \cite{butler_defining_2015}. Note that these days often only return the first day on which the event happens, but do not specify how long the event lasts (in days). 
Rather, these definitions consider days with a certain characteristic that happens consecutively belonging to the same event. 
In these definitions, SSW events are viewed as two distinct events if they are separated by an interval of $n$ days (where the value of $n$ may be 20, 30, or 60, depending on the definition).
For each of these days, we took the corresponding interval from the definitions and found the maximum displacement and split scores in those intervals of days. 
This is a fair comparison: indeed, since we compute the scores daily by day, a day classified as major by some definitions may be just the onset of the event, where the scores are relatively low, even if they ramp up quickly immediately after. 
\begin{figure}[!h]
    \centering
    \includegraphics[scale = 0.38]{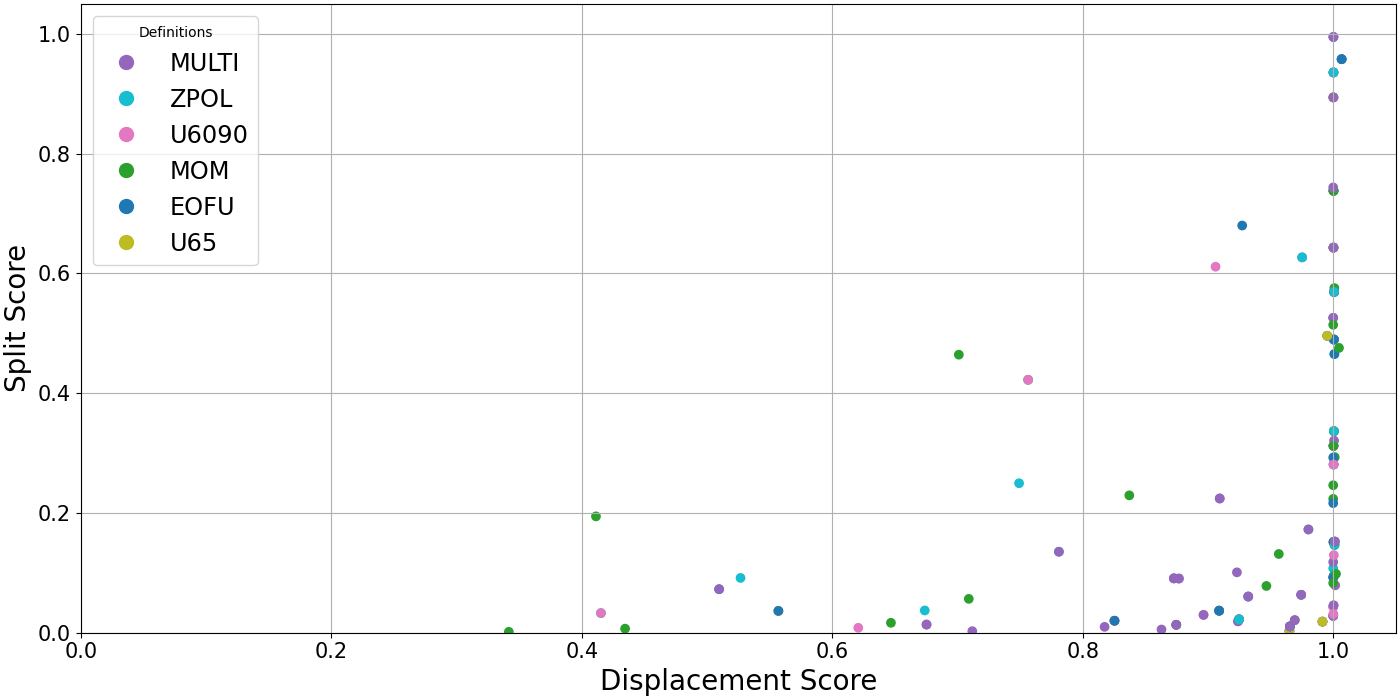}
    \caption{Scatter plot with the displacement (horizontal axis) and split (vertical axis) score, where each point is a collection of days where at least one of the definitions from the literature finds a major event occurring. ``Multi'' stands for days that were classified as major events by multiple definitions.}
    \label{fig:scatter_plot}
\end{figure}
The resulting plot is displayed in Figure~\ref{fig:scatter_plot}. 
For this figure, we observe that scores obtained from our method are consistent with the events identified by previous definitions and shows the corresponding dynamical structures. 
In particular, the average displacement and split scores for the days identified as major by at least one definition in the literature are, respectively, 0.92 and 0.28.
Moreover, most of the events lie in the high displacement score range, and only one day has a score lower than 40\% for one definition, namely the 7th of January, 1986. 
However, when observed, this day is perturbed but never becomes split or displaced, and it is classified as major only by MOM \cite{MOM}, which is quite sensitive to the contour, explaining the classification.
\section{Conclusion and future directions}
\begin{figure}[!h]
    \centering
    \begin{subfigure}{0.49\textwidth}
        \centering
        \includegraphics[scale = 0.39]{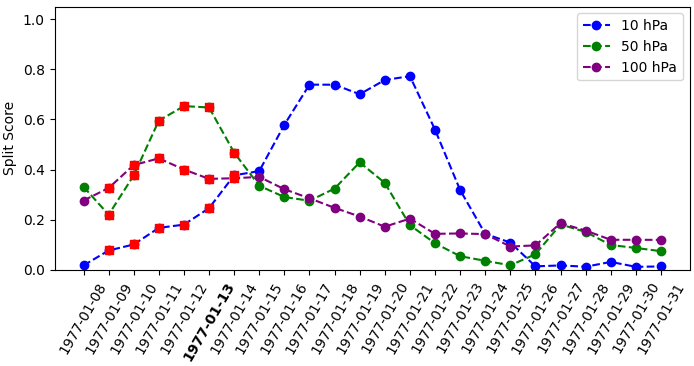}
        \caption{}
    \label{subfig:different_hPa_split}
    \end{subfigure}
     \begin{subfigure}{0.49\textwidth}
        \centering
        \includegraphics[scale = 0.39]{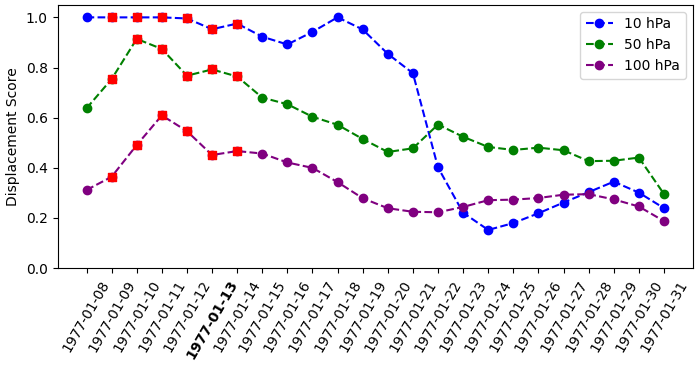}
        \caption{}
    \label{subfig:different_hPa_disp}
    \end{subfigure}
\caption{Split (a) and displacement (b) scores for gph at three different pressures (10,50 and 100 hPa) for the 5 days before and 18 days after 1977-12-13.}
\label{fig:different_hPa}
\end{figure}

In this paper, we propose a definition paradigm for the study of the stratospheric polar vortex, using threshold-free scores computed via persistent homology. 
This novel approach is objective and flexible, allowing for a more consistent picture and the study of the evolution of these events. 
Our code can be used at different pressures in the stratosphere, separately or all at the same time, which opens the way for the study of the vertical structure and anomalies propagation in the stratosphere of these events.
Our method is freely available and very easy to use, and does not require prior knowledge of applied topology. 
Our scores are calculated using the two longest $H_1$-lifespan on the gph scalar field: for how long, scanning the surface downward, do we see a ``circle-like'' hole. 
If there is only one prominent loop or circle-like structure, then the polar vortex is not split, and if this structure is centered away from the North Pole, the vortex is displaced.  Although we do not provide any threshold to define an SSW, users can determine the desired thresholds for weak or major SSWs and displacement or split events based on the SuPerPov scores, as they show the evolution of the polar vortex geometry and strength. (I added a sentence like this here.)
Our approach recovers previous definitions of major SSW events while giving the user a more nuanced picture. The previous definitions use a cut-off threshold to define events, while our approach provides a continuous spectrum and allows the users to inspect the vortex evolution and transition. 
The evolving SuPerPoV scores also provide insights into the prediction of polar vortex states and the related SSWs and, in particular, identify patterns of precursors of the events. 
Moreover, with the prospect of investigating how polar vortex anomalies evolve and propagate in the lower stratosphere and, eventually, on surface weather, our code also includes the possibility of plotting multiple gph for three different pressures at once, namely at 10, 50, and 100 hPa. 
Thus, it can be potentially be used to study the vertical coupling in the stratosphere.
An example is shown in Figure~\ref{fig:different_hPa}, where one can see that, in the split score, there is a peak for the gph at 50hPa right before there is one on the gph at 10 hPa. 
There, we can also see the consistency of the scores evolutions across three different pressure levels.

\section*{Open Research Section}
The data we used is freely available at \cite{our_zenodo}. 
The open source code, together with all information necessary to run it, is available at \cite{github_repo}.

\section*{Conflict of Interest disclosure}
The authors declare no conflicts of interest relevant to this study

\section*{Acknowledgement}
This work was supported by a grant from the Simons Foundation [MPS-TSM-00007525, BG].

\printbibliography

\end{document}